

 \font\smc = cmcsc10
 \voffset = 24pt
 \headline = {{\tenrm\hfil\llap{\number\pageno}}}
 \footline = {\hfil}

\def\section#1{\bigskip\bigskip\bigbreak
 \centerline{{\smc #1}}\nobreak\bigskip\nobreak}
\def\proclaim#1{\medbreak{\smc #1}.\ \begingroup\it}
\def\endproclaim{\endgroup\medbreak}
\def\demo#1{\medbreak{\it #1}.\ }
\def\enddemo{\medbreak}

\def\definition#1{\medbreak{\smc #1}.\ }
\def\enddefinition{\medbreak}
\def\example#1{\medbreak{\smc #1}.\ }
\def\endexample{\medbreak}
\def\title#1{\bigskip\centerline{{\bf{\uppercase{#1}}}}\bigskip}
\def\author#1{\centerline{{\smc #1}}\bigskip}
\def\address#1{\centerline{#1}\bigskip}
\def\email#1{\centerline{E-mail: #1}}
\def\classification#1{\centerline{Classification: #1}}
\def\keywords#1{\centerline{Keywords: #1}}
\def\abstract{\section{Abstract}}
\def\endabstract{\bigbreak}


\def\<{\langle}
\def\>{\rangle}
\def\ADC{{\bf ADC}}
\def\c{\mathbin{\#}}
\def\colim{\mathop{\rm colim}\nolimits}
\def\d{\partial}
\def\epsilon{\varepsilon}
\def\F{{\bf F}}
\def\HOM{\mathop{\rm HOM}\nolimits}
\def\id{\mathop{\rm id}\nolimits}
\def\mod{\ \mathop{\rm mod}\nolimits\ }
\def\N{{\bf N}}
\def\omegacat{\hbox{$\omega$-{\bf cat}}}
\def\Z{{\bf Z}}


\title{Omega-categories and chain complexes}

\author{Richard Steiner}

\address{Department of Mathematics, University of Glasgow,
University Gardens, Glasgow, Scotland G12 8QW}

\email{r.steiner@maths.gla.ac.uk}

\classification{18D05}


\abstract There are several ways to construct omega-categories
from combinatorial objects such as pasting schemes or parity
complexes. We make these constructions into a functor on a
category of chain complexes with additional structure, which we
call augmented directed complexes. This functor from augmented
directed complexes to omega-categories has a left adjoint, and the
adjunction restricts to an equivalence on a category of augmented
directed complexes with good bases. The omega-categories
equivalent to augmented directed complexes with good bases include
the omega-categories associated to globes, simplexes and cubes;
thus the morphisms between these omega-categories are determined
by morphisms between chain complexes. It follows that the entire
theory of omega-categories can be expressed in terms of chain
complexes; in particular we describe the biclosed monoidal
structure on omega-categories and calculate some internal
homomorphism objects.
\endabstract


\section{1. Introduction}

This paper is a contribution to the theory of strict
$\omega$-categories. In the past, $\omega$-categories have been
constructed from combinatorial structures such as pasting schemes,
parity complexes or directed complexes; see Johnson~[6],
Power~[9], Steiner~[10] and Street~[12]. The constructions are not
really functorial, although there are ways to produce non-obvious
morphisms; see Crans-Steiner~[4]. The combinatorial structures can
be regarded as bases for chain complexes, as observed by Kapranov
and Voevodsky in~[7]. We will reinterpret the constructions as a
functor on a category of chain complexes with additional structure
called augmented directed complexes, and we will show that the
functor has a left adjoint. The adjunction is related to the
well-known equivalence between the categories of chain complexes
and of $\omega$-categories in the category of abelian groups; see
Brown-Higgins~[3] for instance.

In the earlier treatments, it is shown that combinatorial
structures which are in a suitable sense loop-free produce free
$\omega$-categories. The functorial version of this result says
that the adjunction restricts to equivalences between certain
pairs of full subcategories. The augmented directed complexes
concerned are free chain complexes with good bases, and the
corresponding $\omega$-categories have good sets of generators.
These $\omega$-categories include those associated to globes,
simplexes and cubes, which determine the entire theory of
$\omega$-categories (see Al-Agl-Brown-Steiner~[1],
Al-Agl-Steiner~[2] and Street~[11]); thus the theory of
$\omega$-categories can be described in terms of chain complexes.
In particular the biclosed monoidal structure on
$\omega$-categories can be described in terms of chain complexes,
and we will calculate some internal homomorphism objects.
Homomorphisms between these $\omega$-categories have been studied
by Gaucher in his work on higher-dimensional automata~[5].

The adjunction between $\omega$-categories and augmented directed
complexes is described in Section~2, bases for augmented directed
complexes are described in Section~3, generating sets for
$\omega$-categories are described in Section~4, the equivalences
between subcategories are described in Section~5, relations with
earlier work are described in Section~6, and the applications to
the theory of $\omega$-categories are described in Section~7.

\section{2. The adjunction}

In this section we describe the adjunction between
$\omega$-categories and augmented directed complexes. First we
define the categories involved. For $\omega$-categories we use the
following description and notation.

\definition{Definition 2.1} An {\it $\omega$-category\/} is a set
with unary source and target operators $d^-_0$, $d^+_0$, $d^-_1$,
$d^+_1$, \dots\ and not everywhere defined binary composition
operators
$$(x,y)\mapsto x\c_0 y,\ (x,y)\mapsto x\c_1 y,\ \ldots$$
such that the following hold:

(i) $x\c_n y$ is defined if and only if $d^+_n x=d^-_n y$;

(ii) for every~$x$ there exists~$n$ such that $d^-_n x=d^+_n x=x$;

(iii) for any~$x$,
$$d^\beta_m d^\alpha_n x=\cases{
 d^\beta_m x& if $m<n$,\cr
 d^\alpha_n x& if $m\geq n$;\cr}$$

(iv) for any~$x$,
$$d^-_n x\c_n x=x\c_n d^+_n x=x;$$

(v) if $x\c_n y$ is defined then
$$\displaylines{
 d^\alpha_m(x\c_n y)
 =d^\alpha_m x
 =d^\alpha_m y
 =d^\alpha_m x\c_n d^\alpha_m y\ {\rm for}\ m<n,\cr
 d^-_n(x\c_n y)=d^-_n x,\cr
 d^+_n(x\c_n y)=d^+_n y,\cr
 d^\alpha_m(x\c_n y)
 =d^\alpha_m x\c_n d^\alpha_m y\ {\rm for}\ m>n;\cr}$$

(vi) for any~$n$,
$$(x\c_n y)\c_n z=x\c_n(y\c_n z)$$
if either side is defined;

(vii) if $m<n$ then
$$(x\c_n y)\c_m(x'\c_n y')=(x\c_m x')\c_n(y\c_m y')$$
when the left side is defined.

A morphism of $\omega$-categories is a function commuting with the
source, target and composition operators. The category of
$\omega$-categories is denoted $\omegacat$.
\enddefinition

Let $C$ be an $\omega$-category. It generates an $\omega$-category
object~$C'$ in the category of abelian groups, which is equivalent
to a nonnegatively graded chain complex~$K$ by~[3]. On~$C'$ there
are two pieces of additional structure: the homomorphism $C'\to\Z$
induced by the morphism from~$C$ to the one-element
$\omega$-category, and the submonoid generated by the image of~$C$
(a submonoid of an abelian group is a subset containing zero and
closed under addition). Correspondingly, it turns out that the
chain complex~$K$ is augmented and has a distinguished submonoid
in each chain group. A chain complex with this kind of additional
structure will be called an augmented directed complex.

\definition{Definition 2.2} An {\it augmented directed
complex\/}~$K$ is an augmented chain complex $(K_n,\d,\epsilon)$
of abelian groups concentrated in nonnegative dimensions, together
with a distinguished submonoid~$K^*_n$ of the chain group~$K_n$
for each~$n$. A morphism of augmented directed complexes from~$K$
to~$L$ is an augmentation-preserving chain map $f\colon K\to L$
such that $f(K^*_n)\subset L^*_n$ for each~$n$. The category of
augmented directed complexes is denoted $\ADC$.
\enddefinition

Similar structures have been studied by Patchkoria in~[8].

In order to construct the functor from $\omegacat$ to $\ADC$, we
first recall the standard filtration for an $\omega$-category.

\proclaim{Proposition 2.3} Let $C$ be an $\omega$-category. Then
$C$~is the union of an increasing sequence of
sub-$\omega$-categories denoted $C_0\subset C_1\subset\ldots\,$,
where
$$C_n
 =d^-_n C
 =d^+_n C
 =\{\,x\in C:d^-_n x=x\,\}
 =\{\,x\in C:d^+_n x=x\,\}.$$
\endproclaim

\demo{Proof} The four definitions for~$C_n$ are consistent because
$d^\beta_n d^\alpha_n=d^\alpha_n$. The~$C_n$ are
sub-$\omega$-categories by Definition 2.1(iii) and~(vi). They form
an increasing sequence because
$d^\alpha_n=d^\beta_{n+1}d^\alpha_n$, and $C$~is the union of the
sub-$\omega$-categories by Definition~2.1(ii).
\enddemo

The filtration of Proposition~2.3 is such that the composition
operators~$\c_m$ are trivial in~$C_n$ for $m\geq n$; thus $C_n$~is
an $n$-category.

We now use the filtration to construct a functor from $\omegacat$
to $\ADC$.

\definition{Definition 2.4} The functor
$\lambda\colon\omegacat\to\ADC$ is defined as follows. Let $C$ be
an $\omega$-category. Then the chain group $(\lambda C)_n$ for
$n\geq 0$ is generated by elements~$[x]_n$ for $x\in C_n$ subject
to relations
$$\hbox{$[x\c_m y]_n=[x]_n+[y]_n$ for $m<n$,}$$
the boundary homomorphism $\d\colon(\lambda C)_{n+1}\to(\lambda
C)_n$ for $n\geq 0$ is given by
$$\d[x]_{n+1}=[d^+_n x]_n-[d^-_n x]_n,$$
the augmentation $\epsilon\colon(\lambda C)_0\to\Z$ is given by
$$\epsilon[x]_0=1,$$
the distinguished submonoid $(\lambda C)_n^*$ is the submonoid
generated by the elements~$[x]_n$.
\enddefinition

To justify this definition, we make the following observations.
First, if $x\in C_{n+1}$ then the $d^\alpha_n x$ are members
of~$C_n$, so the difference $[d^+_n x]_n-[d^-_n x]_n$ is a member
of $(\lambda C)_n$. If $x\c_n y$ is a composite in~$C_{n+1}$ then
$$[d^+_n(x\c_n y)]_n-[d^-_n(x\c_n y)]_n
 =[d^+_n y]_n-[d^-_n x]_n
 =([d^+_n x]_n-[d^-_n x]_n)+([d^+_n y]_n-[d^-_n y]_n)$$
because $d^+_n x=d^-_n y$, and if $x\c_m y$ is a composite
in~$C_{n+1}$ with $m<n$ then
$$\eqalign{
 [d^+_n(x\c_m y)]_n-[d^-_n(x\c_m y)]_n
 &=[d^+_n x\c_m d^+_n y]_n-[d^-_n x\c_m d^-_n y]_n\cr
 &=([d^+_n x]_n+[d^+_n y]_n)-([d^-_n x]_n+[d^-_n y]_n)\cr
 &=([d^+_n x]_n-[d^-_n x]_n)+([d^+_n y]_n-[d^-_n y]_n);\cr}$$
therefore $\d\colon(\lambda C)_{n+1}\to(\lambda C)_n$~is a
well-defined homomorphism. There are no relations on $(\lambda
C)_0$, so $\epsilon\colon(\lambda C)_0\to\Z$ is a well-defined
homomorphism. The composites $\d\d\colon(\lambda
C)_{n+2}\to(\lambda C)_n$ are trivial because $d^\alpha_n
d^-_{n+1}=d^\alpha_n d^+_{n+1}$, and the composite
$\epsilon\d\colon(\lambda C)_1\to\Z$ is obviously trivial.
Finally, $\lambda$~is obviously functorial.

We make the following observation.

\proclaim{Proposition 2.5} If $C$~is an $\omega$-category and
$x\in C_m$ with $m<n$, then $[x]_n=0$ in $(\lambda C)_n$.
\endproclaim

\demo{Proof} This holds because
$$[x]_n=[x\c_m d^+_m x]_n=[x\c_m x]_n=[x]_n+[x]_n.$$
\enddemo

We will now go from chain complexes to $\omega$-categories. We
first define a functor~$\mu$ from arbitrary chain complexes to
$\omega$-categories in the category of abelian groups. It is the
equivalence given implicitly in~[3], and it may also be regarded
an additive version of Street's construction~[12].

\definition{Definition 2.6} The functor~$\mu$ from chain complexes
to $\omega$-categories in the category of abelian groups is
defined as follows. Let $K$ be a chain complex. Then $\mu K$ is
the abelian group of double sequences
$$x=(x^-_0,x^+_0,x^-_1,x^+_1,\ldots\,)$$
such that
$$\eqalign{
 &\hbox{$x^-_n\in K_n$ and $x^+_n\in K_n$,}\cr
 &\hbox{$x^-_n=x^+_n=0$ for all but finitely many values
 of~$n$,}\cr
 &\hbox{$x^+_n-x^-_n=\d x^-_{n+1}=\d x^+_{n+1}$ for
 $n\geq 0$;}\cr}$$
if $x\in\mu K$ then
$$d^\alpha_n x
 =(x^-_0,x^+_0,\ldots,
 x^-_{n-1},x^+_{n-1},x^\alpha_n,x^\alpha_n,0,0,\ldots\,);$$
if $d^+_n x=d^-_n y=z$, say, in $\mu K$ then
$$\eqalign{
 x\c_n y
 &=x-z+y\cr
 &=(x^-_0,\ y^+_0,\ \ldots,\ x^-_n,\ y^+_n,\
 x^-_{n+1}+y^-_{n+1},\ x^+_{n+1}+y^+_{n+1},\ \ldots\,).\cr}$$
\enddefinition

It is straightforward to check that $\mu$~is a well-defined
functor. We extend the filtration of $\mu K$ given in
Proposition~2.3 by writing $(\mu K)_{-1}=0$, and we get the
following result.

\proclaim{Proposition 2.7} Let $K$ be a chain complex and let $n$
be a nonnegative integer. Then $x^-_n=x^+_n$ for $x\in (\mu K)_n$,
and the homomorphism $x\mapsto x^\alpha_n\colon(\mu K)_n\to K_n$
fits into a natural split short exact sequence
$$0\to(\mu K)_{n-1}\to(\mu K)_n\to K_n\to 0.$$
\endproclaim

\demo{Proof} If $x\in(\mu K)_n$ then $x=d^\alpha_n x$, so
$x^-_n=x^+_n$. For $x\in\mu K$, one finds that $x\in(\mu K)_n$ if
and only if $x^-_m=x^+_m=0$ for all $m>n$, and it follows that the
homomorphism $x\mapsto x^\alpha_n\colon(\mu K)_n\to K_n$ has
kernel $(\mu K)_{n-1}$. This means that there is a natural exact
sequence $0\to(\mu K)_{n-1}\to(\mu K)_n\to K_n$. To complete the
proof we need a natural splitting homomorphism $K_n\to(\mu K)_n$.
For $n=0$ we use the homomorphism
$$x\mapsto(x,x,0,0,\ldots\,),$$
and for $n>0$ we can use the homomorphism
$$x\mapsto(0,0,\ldots,0,\d x,x,x,0,0,\ldots\,).$$
\enddemo

We now use the additional structure on an augmented directed
complex to define a functor from $\ADC$ to $\omegacat$.

\definition{Definition 2.8} The functor
$\nu\colon\ADC\to\omegacat$ is defined as follows. Let $K$ be an
augmented directed complex. Then $\nu K$ is the
sub-$\omega$-category of $\mu K$ consisting of the elements
$$(x^-_0,x^+_0,x^-_1,x^+_1,\ldots\,)$$
such that
$$\eqalign{
 &\hbox{$x^-_n\in K^*_n$ and $x^+_n\in K^*_n$ for all~$n$,}\cr
 &\epsilon x^-_0=\epsilon x^+_0=1.\cr}$$
\enddefinition

It is straightforward to check that Definition~2.8 gives a
well-defined functor from augmented directed complexes to
$\omega$-categories: if $K$~is an augmented directed complex then
$\nu K$ is a subset of $\mu K$ closed under the $\omega$-category
operations. Obviously $\nu K$ is not a subgroup of $\mu K$.

Finally in this section, we show that $\lambda$~is left adjoint
to~$\nu$. First we describe the unit of the adjunction.

\definition{Definition 2.9} The natural transformation
$\eta\colon C\to\nu\lambda C$ for an $\omega$-category~$C$ is
defined by the formula
$$\eta x
 =([d^-_0 x]_0,\ [d^+_0 x]_0,\ [d^-_1 x]_1,\ [d^+_1 x]_1,\
 \ldots\,).$$
\enddefinition

In order to justify this definition, we first show that $\eta
x\in\nu\lambda C$ for $x\in C$. Indeed, it is clear that
$[d^\alpha_n x]_n\in(\lambda C)^*_n$ and $\epsilon[d^\alpha_0
x]_0=1$, and we also have
$$\d[d^\alpha_{n+1}x]_{n+1}
 =[d^+_n d^\alpha_{n+1}x]_n-[d^-_n d^\alpha_{n+1}x]_n
 =[d^+_n x]_n-[d^-_n x]_n,$$
so it suffices to show that $[d^\alpha_n x]_n=0$ for $n$
sufficiently large. But by Proposition~2.3 there exists~$p$ such
that $d^\alpha_n x\in C_p$ for all~$n$, and it then follows from
Proposition~2.5 that $[d^\alpha_n x]_n=0$ for $n>p$.

We must also show that $\eta\colon C\to\nu\lambda C$ is a morphism
of $\omega$-categories. But we get $\eta d^\alpha_n
x=d^\alpha_n\eta x$ because $[d^\beta_m d^\alpha_n x]_m=[d^\beta_m
x]_m$ for $m<n$, because $[d^\beta_n d^\alpha_n x]_n=[d^\alpha_n
x]_n$, and because $[d^\beta_m d^\alpha_n x]_m=[d^\alpha_n x]_m=0$
for $m>n$ by Proposition~2.5. We also get $\eta(x\c_n y)=\eta
x\c_n\eta y$ because $[d^-_m(x\c_n y)]_m=[d^-_m x]_m$ and
$[d^+_m(x\c_n y)]_m=[d^+_m y]_m$ for $m\leq n$ and because
$[d^\alpha_m(x\c_n y)]_m=[d^\alpha_m x\c_n d^\alpha_m y]_m
=[d^\alpha_m x]_m+[d^\alpha_m y]_m$ for $m>n$. It follows that
$\eta$~is a morphism of $\omega$-categories, and it is clear that
$\eta$~is natural.

Next we describe the counit. We denote this by~$\pi$ rather
than~$\epsilon$, in order to avoid confusion with augmentations.

\definition{Definition 2.10} The natural transformation
$\pi\colon\lambda\nu K\to K$ is defined for an augmented directed
complex~$K$ by the formula
$$\pi[x]_n=x^-_n=x^+_n.$$
\enddefinition

To justify this definition, observe first that $x^-_n=x^+_n$ for
$x\in(\nu K)_n$ by Proposition~2.7. The formula
$\pi[x]_n=x^-_n=x^+_n$ gives a well-defined homomorphism
$\pi\colon(\lambda\nu K)_n\to K_n$ because if $m<n$ then $(x\c_m
y)^\alpha_n=x^\alpha_n+y^\alpha_n$. We get a chain map because
$$\d\pi[x]_{n+1}
 =\d x^\alpha_{n+1}
 =x^+_n-x^-_n
 =(d^+_n x)^+_n-(d^-_n x)^-_n
 =\pi([d^+_n x]_n-[d^-_n x]_n)
 =\pi\d[x]_{n+1}.$$
This chain map is augmentation-preserving because
$\epsilon\pi[x]_0=\epsilon x^\alpha_0=1=\epsilon[x]_0$ for
$x\in(\nu K)_0$, and we get $\pi(\lambda\nu K)^*_n\subset K^*_n$
because $x^\alpha_n\in K^*_n$ for $x\in\nu K$. Therefore
$\pi\colon\lambda\nu K\to K$ is a morphism of augmented directed
complexes. It is clearly natural.

The main theorem is now as follows.

\proclaim{Theorem 2.11} The functors
$\lambda\colon\omegacat\to\ADC$ and $\nu\colon\ADC\to\omegacat$
form an adjoint pair with unit $\eta\colon\id\to\nu\lambda$ and
counit $\pi\colon\lambda\nu\to\id$.
\endproclaim

\demo{Proof} We must show that
$$(\pi\lambda)\circ(\lambda\eta)=\id\colon\lambda\to\lambda$$
and
$$(\nu\pi)\circ(\eta\nu)=\id\colon\nu\to\nu.$$

Let $C$ be an $\omega$-category. The generators of $(\lambda C)_n$
have the form~$[x]_n$ with $x\in C_n$. For these generators,
$d^\alpha_n x=x$, so
$$\eqalign{
 \pi(\lambda\eta)[x]_n
 &=\pi[\eta x]_n\cr
 &=\pi\bigl[([d^-_0 x]_0,\ [d^+_0 x]_0,\ [d^-_1 x]_1,\ [d^+_1
 x]_1,\ \ldots\,)\bigr]_n\cr
 &=[d^\alpha_n x]_n\cr
 &=[x]_n;\cr}$$
therefore $(\pi\lambda)\circ(\lambda\eta)=\id$.

Now let $K$ be an augmented directed complex and let $x$ be a
member of $\nu K$. Then
$$\eqalign{
 (\nu\pi)\eta x
 &=(\nu\pi)([d^-_0 x]_0,\ [d^+_0 x]_0,\ [d^-_1 x]_1,\ [d^+_1
 x]_1,\ \ldots\,)\cr
 &=(\pi[d^-_0 x]_0,\ \pi[d^+_0 x]_0,\ \pi[d^-_1 x]_1,\ \pi[d^+_1
 x]_1,\ \ldots\,)\cr
 &=\bigl((d^-_0 x)^-_0,\ (d^+_0 x)^+_0,\ (d^-_1 x)^-_1,\ (d^+_1
 x)^+_1,\ \ldots\,\bigr)\cr
 &=(x^-_0,x^+_0,x^-_1,x^+_1,\ldots\,)\cr
 &=x;\cr}$$
therefore $(\nu\pi)\circ(\eta\nu)=\id$.

This completes the proof.
\enddemo

\section{3. Bases for augmented directed complexes}

We will now consider augmented directed complexes with bases. We
essentially recover the examples constructed in earlier treatments
such as [6], [9], [10] and~[12].

\definition{Definition 3.1} Let $K$ be an augmented directed
complex. A {\it basis\/} for~$K$ is a set $B\subset\bigsqcup_n
K_n$ such that each~$K_n$ is a free abelian group with basis
$B\cap K_n$ and each~$K^*_n$ is the submonoid of~$K_n$ generated
by $B\cap K_n$.
\enddefinition

Suppose that $K$~is an augmented directed complex with a basis. We
make~$K_n$ into a partially ordered abelian group by the rule
$$x\leq y\iff y-x\in K^*_n.$$
The basis elements in~$K_n$ can be characterised as the minimal
non-zero elements in~$K^*_n$, and it follows that $K$~has only one
basis. Note also that $K_n$~is a lattice: any two elements
$x$~and~$y$ have a least upper bound $x\vee y$ and a greatest
lower bound $x\wedge y$. If $x$~is an element of~$K_{n+1}$, then
there are unique elements $\d^- x,\d^+ x\in K_n$, the negative and
positive parts of $\d x$, such that
$$\d x=\d^+ x-\d^- x,\quad \d^- x\wedge\d^+ x=0;$$
indeed, if $\d x$ is expressed as a linear combination of distinct
basis elements, then $\d^+ x$ is the sum of the terms with
positive coefficients and $-\d^- x$ is the sum of the terms with
negative coefficients.

Let $b$ be a basis element for~$K$. We denote the {\it
dimension\/} of~$b$ by~$|b|$, so that $b\in K_{|b|}$. We then
define elements $\<b\>^-_n$~and~$\<b\>^+_n$ in~$K_n$ by downward
recursion as follows:
$$\<b\>^\alpha_n=\cases{
 0& for $n>|b|$,\cr
 b& for $n=|b|$,\cr
 \d^\alpha\<b\>^\alpha_{n+1}& for $n<|b|$.\cr}$$
It is straightforward to check that this produces an
element~$\<b\>$ of $\mu K$; in fact we can make the following
definition.

\definition{Definition 3.2} Let $K$ be an augmented directed
complex with a basis and let $b$ be a basis element. Then the {\it
atom\/} associated to~$b$ is the element~$\<b\>$ of $\mu K$ such
that $\<b\>^\alpha_n=0$ for $n>|b|$, such that
$\<b\>^\alpha_{|b|}=b$, and such that $\<b\>^-_n\wedge\<b\>^+_n=0$
for $n<|b|$. The {\it dimension\/} of the atom~$\<b\>$ is the
dimension of the corresponding basis element~$b$.
\enddefinition

From Proposition~2.7 we deduce the following result.

\proclaim{Proposition 3.3} Let $K$ be an augmented directed
complex with a basis. Then the atoms form a basis for the abelian
group $\mu K$, the $n$-dimensional atoms form a basis for $(\mu
K)_n/(\mu K)_{n-1}$, and the atoms of dimension greater than~$n$
form a basis for $\mu K/(\mu K)_n$.
\endproclaim

If $\<b\>$~is an atom, then $\<b\>^\alpha_n\in K_n^*$ for all~$n$
by construction. An atom~$\<b\>$ is therefore in $\nu K$ if and
only if $\epsilon\<b\>^-_0=\epsilon\<b\>^+_0=1$. This leads us to
the following definition.

\definition{Definition 3.4} A basis~$B$ for an augmented directed
complex is {\it unital\/} if
$\epsilon\<b\>^-_0=\epsilon\<b\>^+_0=1$ for every $b\in B$.
\enddefinition

For the equivalence theorem of Section~5 we need bases which are
unital and are also loop-free in the sense of the following
definition.

\definition{Definition 3.5} A basis~$B$ for an augmented directed
complex is {\it loop-free\/} if there are partial orderings
$\leq_0,\,\leq_1,\,\ldots\,$ on~$B$ such that $a<_n b$ whenever
$\<a\>^+_n\wedge\<b\>^-_n>0$ and $|a|,|b|>n$.
\enddefinition

Note that the partial orderings~$\leq_n$ are quite different from
the partial orderings on the individual chain groups.

In practice, one usually has a stronger condition.

\definition{Definition 3.6} A basis~$B$ for an augmented directed
complex is {\it strongly loop-free\/} if there is a partial
ordering~$\leq_\N$ on~$B$ such that $a<_\N b$ whenever
$a\leq\d^-b$ or $\d^+a\geq b$.
\enddefinition

The two notions of loop-freeness are related as follows.

\proclaim{Proposition 3.7} If a basis for an augmented directed
complex is strongly loop-free, then it is loop-free.
\endproclaim

\demo{Proof} Let $\leq_\N$ be a partial ordering with the property
required for strong loop-freeness. For each~$n$, we will show that
$\leq_\N$~has the property required for~$\leq_n$ in the definition
of loop-freeness. In other words, we suppose that
$\<a\>^+_n\wedge\<b\>^-_n>0$ with $|a|,|b|>n$ for some~$n$, and we
show that $a<_\N b$. Indeed we can choose a basis element~$c$ with
$\<a\>^+_n\geq c$ and $c\leq\<b\>^-_n$, and we will show that
$a<_\N c<_\N b$.

To show that $c<_\N b$, observe that, by the construction
of~$\<b\>$, we must have $c\leq\d^-c'$ for some basis element~$c'$
with $c'\leq\<b\>^-_{n+1}$, and we than have $c<_\N c'$. If
$n+1<|b|$ then we repeat this argument, and eventually we get
$c<_\N\ldots<_\N c''$ with $c''\leq\<b\>^-_{|b|}$. But
$\<b\>^-_{|b|}=b$, so $c''=b$, and we have got $c<_\N b$ as
claimed. The proof that $a<_\N c$ is similar.
\enddemo

\example{Example 3.8} Let $K$ be the chain complex of a simplicial
set. Then $K$~is an augmented chain complex with a distinguished
basis, so it can be regarded as an augmented directed complex with
a basis. In particular, let $\Delta[p]$ be the chain complex of
the standard $p$-simplex, so that $\Delta[p]$ has the following
structure: the basis elements are the ordered $(n+1)$-tuples of
integers $(v_0,\ldots,v_n)$ with $0\leq v_0<v_1<\ldots<v_n\leq p$;
the dimension of $(v_0,\ldots,v_n)$ is~$n$; the boundary
$\d\colon\Delta[p]_n\to\Delta[p]_{n-1}$ is given by
$$\d=\d_0-\d_1+\d_2-\ldots+(-1)^n\d_n,$$
where
$$\d_i(v_0,\ldots,v_n)=(v_0,\ldots,v_{i-1},v_{i+1},\ldots,v_n);$$
the augmentation is given by $\epsilon(v_0)=1$. For $m\leq n$ one
finds that
$$\<(v_0,\ldots,v_n)\>^\alpha_m
 =\sum\d_{i(1)}\ldots\d_{i(n-m)}(v_0,\ldots,v_n),$$
where the sum runs over $(n-m)$-tuples such that $0\leq
i(1)<\ldots<i(n-m)\leq n$ and the parities
$(-)^{i(1)},\ldots,(-)^{i(n-m)}$ form the alternating sequence
$\alpha,-\alpha,\alpha,\ldots\,$. We have therefore recovered
Street's oriented simplexes~[11]. The basis is unital, because
$\<(v_0,\ldots,v_n)\>^-_0=(v_0)$ and
$\<(v_0,\ldots,v_n)\>^+_0=(v_n)$. The basis is also strongly
loop-free under the total ordering given recursively as follows:
$(v_0,\ldots,v_n)<_\N(w_0,\ldots,w_m)$ if
$$v_0<w_0,$$
or if
$$\hbox{$v_0=w_0$ and $n=0$ and $m>0$,}$$
or if
$$\hbox{$v_0=w_0$ and $n>0$ and $m>0$ and
 $(v_1,\ldots,v_n)>_\N(w_1,\ldots,w_m)$.}$$
\endexample

\example{Example 3.9} Let $B$ be a finite non-empty totally
ordered set in which each element~$b$ is assigned a nonnegative
integer dimension~$|b|$. Suppose also that the initial and final
elements have dimension~$0$ and that adjacent elements have
dimensions differing by~$1$. For $b\in B$ with $|b|>0$, let
$\delta^- b$ be the last element of dimension $|b|-1$ to come
before~$b$, and let $\delta^+ b$ be the first element of dimension
$|b|-1$ to come after~$b$; the hypotheses ensure that these
elements exist. The hypotheses also ensure that the elements
between $\delta^- b$ and $\delta^+ b$ have dimension at
least~$|b|$, and for $|b|>1$ it follows that
$\delta^\alpha\delta^- b=\delta^\alpha\delta^+ b$. By taking $\d
b=\delta^+ b-\delta^- b$ for $|b|>0$ and $\epsilon b=1$ for
$|b|=0$ we define an augmented directed complex with basis~$B$.
Clearly $\<b\>^\alpha_n$ is a basis element for $n\leq |b|$, so
the basis is unital. It is also strongly loop-free under the
original total ordering.

In particular, let $p$ be a nonnegative integer and let the
sequence of dimensions of the elements of~$B$ be
$$0,1,\ldots,p-1,p,p-1,\ldots,1,0.$$
Then the augmented directed complex is called the {$p$-dimensional
globe\/} and denoted $G[p]$. Let $x$ denote the $p$-dimensional
basis element; then $d^-_p\<x\>=d^+_p\<x\>=\<x\>$, and the atoms
other than~$\<x\>$ are the elements $d^\alpha_i\<x\>$ for $i<p$.

A $p$-dimensional globe is in a sense free on a $p$-dimensional
generator. In a similar way we can get an augmented directed
complex $G[p;n]$ free on a $\c_n$-composable pair of
$p$-dimensional elements: we take the sequence of dimensions of
the elements of~$B$ to run from~$0$ up to~$p$, then down to
$\min\{p,n\}$, then up to~$p$, then down to~$0$. Let the
$p$-dimensional basis elements in order be $x$~and~$y$ (they
coincide if $p\leq n$); then $d^-_p\<x\>=d^+_p\<x\>=\<x\>$,
$d^-_p\<y\>=d^+_p\<y\>=\<y\>$, and $d^+_n\<x\>=d^-_n\<y\>$. The
atoms other than $\<x\>$~and~$\<y\>$ are the $d^\alpha_i\<x\>$ and
$d^\alpha_i\<y\>$ for $i<p$, with the identifications
$d^\alpha_i\<x\>=d^\alpha_i\<y\>$ for $i<\min\{p,n\}$ and with the
identification $d^+_n\<x\>=d^-_n\<y\>$ if $n<p$.

There is a similar augmented directed complex $G[p;n,n]$ free on a
$\c_n$-composable triple of $p$-dimensional elements: the sequence
of dimensions is now from~$0$ up to~$p$, down to $\min\{p,n\}$, up
to~$p$, down to $\min\{p,n\}$, up to~$p$, down to~$0$. For $m<n$
there is also an augmented directed complex $G[p;n,m,n]$ free on a
composable quadruple of $p$-dimensional elements in the
configuration $(x\c_n y)\c_m(x'\c_n y')$ of Definition~2.1(vii):
the sequence of dimensions is now from~$0$ up to~$p$, down to
$\min\{p,n\}$, up to~$p$, down to $\min\{p,m\}$, up to~$p$, down
to $\min\{p,n\}$, up to~$p$, down to~$0$.
\endexample

\example{Example 3.10} The category of augmented chain complexes
has a symmetric monoidal structure under the tensor product
$(K,L)\mapsto K\otimes L$, where
$$\eqalign{
 &(K\otimes L)_n=\bigoplus_i K_i\otimes L_{n-i},\cr
 &\d(x\otimes y)=\d x\otimes y+(-1)^{|x|}x\otimes\d y,\cr
 &\epsilon(x\otimes y)=(\epsilon x)(\epsilon y)\ {\rm for\ }
 |x|=|y|=0\cr}$$
(we write $|x|=i$ if $x\in K_i$, etc.). The identity object is the
$0$-dimensional globe $G[0]$ of Example~3.9. We extend this
structure to a monoidal structure on augmented directed complexes
as follows: $(K\otimes L)_n^*$ is the submonoid of $(K\otimes
L)_n$ generated by the elements $x\otimes y$ with $x\in K_i^*$ and
$y\in L_{n-i}^*$. Note that the monoidal structure on augmented
directed complexes is not symmetric: the standard switch morphism
$x\otimes y\mapsto(-1)^{|x|\,|y|}(y\otimes x)$ does not always map
$(K\otimes L)_n^*$ into $(L\otimes K)_n^*$.

Suppose now that $K$~and~$L$ are augmented directed complexes with
bases $A$~and~$B$. Then $K\otimes L$ has a basis~$C$ consisting of
the elements $a\otimes b$ for $a\in A$ and $b\in B$. One finds
that
$$\<a\otimes b\>^\alpha_n
 =\sum_{i=0}^n\<a\>_i^\alpha\otimes\<b\>^{(-)^i\alpha}_{n-i};$$
in particular $\<a\otimes
b\>^\alpha_0=\<a\>^\alpha_0\otimes\<b\>^\alpha_0$. If $A$~and~$B$
are unital, then it follows that $C$~is unital. Similarly, if
$A$~and~$B$ are strongly loop-free under partial
orderings~$\leq_\N$, then $C$~is strongly loop-free under the
partial ordering~$\leq_\N$ such that $a\otimes b\leq_\N a'\otimes
b'$ for
$$a<_\N a',$$
or for
$$\hbox{$a=a'$ and $|a|$ even and $b\leq_\N b'$,}$$
or for
$$\hbox{$a=a'$ and $|a|$ odd and $b\geq_\N b'$.}$$

In particular, for $p\geq 0$, let $Q[p]$ be the chain complex of
the $p$-dimensional cube. Then $Q[p]$ is the $p$-fold tensor power
of the one-dimensional globe $G[1]$, so $Q[p]$ is an augmented
directed complex with a strongly loop-free unital basis.
\endexample

\section{4. Bases for $\omega$-categories}

We will now describe bases for $\omega$-categories, corresponding
to bases for augmented directed complexes. The required properties
can be described directly in terms of $\omega$-categories, but it
is often easier to work in the augmented directed complexes got by
applying~$\lambda$.

We begin with a particular kind of generating set, analogous to a
spanning set. Recall that if $C$~is an $\omega$-category then
$C_n$~is the sub-$\omega$-category $d^-_n C=d^+_n C$.

\definition{Definition 4.1} An $\omega$-category~$C$ is
{\it composition-generated\/} by a subset~$E$ if each member~$e$
of~$E$ is assigned a nonnegative integer dimension~$|e|$ and if
for $n\geq 0$ the sub-$\omega$-category~$C_n$ is generated under
the composition operations~$\c_m$ by the elements of~$E$ of
dimension at most~$n$.
\enddefinition

Composition-generation should be distinguished from ordinary
generation, where one uses the operations~$d^\alpha_m$ as well as
the composition operations.

Note that if $e$~is a composition-generator for an
$\omega$-category~$C$ then $e\in C_{|e|}$, so that $d^\alpha_n
e=e$ for $n\geq |e|$.

There are standard forms for the elements of $\omega$-categories
with composition-generators, as follows.

\proclaim{Proposition 4.2} Let $C$ be an $\omega$-category with a
set of composition-generators and let $x$ be a member of~$C$. Then
$x$~is a generator or $x$~has an expression $x=x_1\c_r\ldots\c_r
x_k$ with $r\geq 0$ and $k\geq 2$ such that the~$x_i$ are
composites of generators, each~$x_i$ has exactly one factor of
dimension greater than~$r$, and at most one of the~$x_i$ has a
factor of dimension greater than $r+1$.
\endproclaim

\demo{Proof} It suffices to prove the result when $x$~is a
composite, say $x=y\c_n z$, in which $y$~and~$z$ have expressions
of the required form. If the expression for~$y$ involves no
factors of dimension greater than~$n$, then $y$~is an identity
for~$\c_n$ by Definition~2.1(v), so $x=z$ and the result therefore
holds for~$x$. A similar argument applies if the expression
for~$z$ involves no factors of dimension greater than~$n$.

From now on, assume that the expressions for $y$~and~$z$ both
involve factors of dimension greater than~$n$. Let $r$ be the
largest integer such that the expressions for $y$~and~$z$ together
involve more than one factor of dimension greater than~$r$; thus
$r\geq n$. Suppose that the expressions for $y$~and~$z$ have
$k$~and~$l$ factors respectively of dimension greater than~$r$;
thus $k+l\geq 2$. Then there is a decomposition
$$y=y_1\c_r\ldots\c_r y_k\c_r d^+_r y\c_r\ldots\c_r d^+_r y,$$
with $l$~appearances of $d^+_r y$, such that $y_i$~has exactly one
factor of dimension greater than~$r$ (if $k=0$ then $y=d^+_r y$,
so this holds trivially; if $k=1$ we take $y_1=y$; if $k>1$ then
$y_1\c_r\ldots\c_r y_k$ is the given expression for~$y$).
Similarly there is a decomposition
$$z=d^-_r z\c_r\ldots\c_r d^-_r z\c_r z_1\c_r\ldots\c_r z_l,$$
with $k$~factors equal to $d^-_r z$, such that $z_j$~has exactly
one factor of dimension greater than~$r$. The choice of~$r$
ensures that the $y_i$~and~$z_j$ have at most one factor of
dimension $r+1$ between them. If now $n=r$ then the desired
decomposition of $x=y\c_r z$ is given by
$$x=y_1\c_r\ldots\c_r y_k\c_r z_1\c_r\ldots\c_r z_l;$$
if $n<r$ then Definition~2.1(vii) gives
$$\eqalign{
 x
 &=(y_1\c_r\ldots\c_r y_k\c_r d^+_r y\c_r\ldots\c_r d^+_r y)
 \c_n
 (d^-_r z\c_r\ldots\c_r d^-_r z\c_r z_1\c_r\ldots\c_r z_l)\cr
 &=(y_1\c_n d^-_r z)\c_r\ldots\c_r(y_k\c_n d^-_r z)
 \c_r(d^+_r y\c_n z_1)\c_r\ldots\c_r(d^+_r y\c_n z_l),\cr}$$
which is a decomposition of the required form because $d^-_r z$
and $d^+_r y$ are composites of generators of dimension at
most~$r$.

This completes the proof.
\enddemo

Composition-generators for an $\omega$-category~$C$ produce
generators for $(\lambda C)_n$ and $(\lambda C)_n^*$ as follows.

\proclaim{Proposition 4.3} Let $C$ be an $\omega$-category with a
set of composition-generators~$E$. Then the abelian group
$(\lambda C)_n$ and the submonoid $(\lambda C)_n^*$ are generated
by the elements~$[e]_n$ for $e$ an $n$-dimensional element of~$E$.
\endproclaim

\demo{Proof} The abelian group $(\lambda C)_n$ and the submonoid
$(\lambda C)_n^*$ are generated by the elements~$[x]_n$ for $x\in
C_n$, and the elements of~$C_n$ are the composites of the members
of~$E$ of dimension at most~$n$. If $m\geq n$ then the elements
of~$C_n$ are identities for~$\c_m$, because $d^\beta_m
d^\alpha_n=d^\alpha_n$; it therefore suffices to use the
composition operators~$\c_m$ for $m<n$. These operations become
addition in $(\lambda C)_n$, so $(\lambda C)_n$ and $(\lambda
C)_n^*$ are generated by the elements~$[e]_n$ for $e\in E$ with
$|e|\leq n$. Finally $[e]_n=0$ for $|e|<n$ by Proposition~2.5;
therefore $(\lambda C)_n$ and $(\lambda C)_n^*$ are generated by
the elements~$[e]_n$ for $e\in E$ with $|e|=n$.
\enddemo

We now define a basis for an $\omega$-category~$C$ in terms of a
basis for the augmented directed complex $\lambda C$.

\definition{Definition 4.4} A {\it basis\/} for an
$\omega$-category~$C$ is a set of composition-generators~$E$ such
$\lambda C$ has a basis and the function $e\mapsto[e]_{|e|}$
maps~$E$ bijectively onto the basis for $\lambda C$.
\enddefinition

Because of Proposition~4.3, a composition-generating set~$E$ is a
basis if and only if the elements~$[e]_{|e|}$ for $e\in E$ are
distinct and linearly independent.

Let $C$ be an $\omega$-category with a basis~$E$. Then $\lambda C$
has a basis~$B$, and we can use properties of~$B$ to define
properties of~$E$.

\definition{Definition 4.5} Let $E$ be a basis for an
$\omega$-category~$C$. Then

(i) $E$~is {\it atomic\/} if $[d^-_n e]_n\wedge[d^+_n e]_n=0$ for
$e\in E$ and $n<|e|$;

(ii) $E$~is {\it loop-free\/} if the basis for $\lambda C$ is
loop-free;

(iii) $E$~is {\it strongly loop-free\/} if the basis for $\lambda
C$ is strongly loop-free.
\enddefinition

The point of atomicity is as follows.

\proclaim{Proposition 4.6} Let $E$ be an atomic basis for an
$\omega$-category~$C$. Then $\eta e=\<[e]_{|e|}\>$ for each~$e$
in~$E$, and the basis for $\lambda C$ is unital.
\endproclaim

\demo{Proof} Recall from Definition~2.9 that $(\eta
e)^\alpha_n=[d^\alpha_n e]_n$ for all $n$~and~$\alpha$. If $n>|e|$
then $[d^\alpha_n e]_n=[e]_n=0$ by Proposition~2.5. If $n=|e|$
then $[d^\alpha_n e]_n=[e]_{|e|}$. By comparing Definition~4.5(i)
with Definition~3.2, we now see that $\eta e=\<[e]_{|e|}\>$. Since
$\eta e\in\nu\lambda C$, it follows that
$\epsilon\<[e]_{|e|}\>^\alpha_0=\epsilon(\eta e)^\alpha_0=1$;
therefore the basis for $\lambda C$ is unital. This completes the
proof.
\enddemo

\example{Example 4.7} Let $p$ be a nonnegative integer, and let
$F[p]$ be the $\omega$-category with the following presentation:
there is a single generator~$u$ and there are relations $d^-_p
u=d^+_p u=u$. The $\omega$-categories $F[p]$ represent the
elements of $\omega$-categories, in the sense that there are
natural bijections between $\hom(F[p],C)$ and~$C_p$. We have
$d^\alpha_i u=u$ for $i\geq p$. Using Definition 2.1(iii) and~(v),
we see that $F[p]$ is composition-generated by~$u$ and the
elements $d^\alpha_i u$ for $i<p$.

Now let $G[p]$ be the $p$-dimensional globe with $p$-dimensional
basis element~$x$ as in Example~3.9. The atom~$\<x\>$ in $\nu
G[p]$ satisfies the relations $d^-_p\<x\>=d^+_p\<x\>=\<x\>$, so
there is a morphism $F[p]\to\nu G[p]$ given by $u\mapsto\<x\>$.
The adjoint $\lambda F[p]\to G[p]$ maps the generators $[u]_p$ and
$[d^\alpha_i u]_i$ for $\lambda F[p]$ (see Proposition~4.3)
bijectively to the basis elements for $G[p]$, and it follows that
the morphism $\lambda F[p]\to G[p]$ is an isomorphism. It then
follows that the elements $u$ and $d^\alpha_i u$ form a basis for
$F[p]$, and one can check that this basis is atomic. It is also
strongly loop-free, because the basis for $G[p]$ is strongly
loop-free.

There are similar $\omega$-categories $F[p;n]$, $F[p;n,n]$ and
$F[p;n,m,n]$ (where $m<n$) corresponding to the other augmented
directed complexes of Example~3.9. These $\omega$-categories all
have strongly loop-free atomic bases. There is a presentation for
$F[p;n]$ given by
$$\<\,u,v: d^\alpha_p u=u,\ d^\alpha_p v=v,\ d^+_n u=d^-_n v\,\>,$$
there is a presentation for $F[p;n,n]$ given by
$$\<\,u,v,w: d^\alpha_p u=u,\ d^\alpha_p v=v,\ d^\alpha_p w=w,\
 d^+_n u=d^-_n v,\ d^+_n v=d^-_n w\,\>,$$
and there is a presentation for $F[p;n,m,n]$ given by
$$\<\,u,v,u',v':d^\alpha_p u=u,\ d^\alpha_p v=v,\
 d^\alpha_p u'=u',\ d^\alpha_p v'=v',\
 d^+_n u=d^-_n v,\ d^+_m v=d^-_m u',\ d^+_n u'=d^-_n v'\,\>.$$
\endexample

\section{5. The adjoint equivalence}

From Definition~4.4, Definition~4.5 and Proposition~4.6, if an
$\omega$-category~$C$ has a loop-free atomic basis then $\lambda
C$ has a loop-free unital basis. We will now show that
$\omega$-categories with loop-free atomic bases and augmented
directed complexes with loop-free unital bases are equivalent
under the adjoint functors $\lambda$~and~$\nu$.

We begin by considering decompositions in $\nu K$, where $K$~is an
augmented directed complex with a basis. The basic existence
result is as follows.

\proclaim{Proposition 5.1} Let $K$ be an augmented directed
complex with a basis and let $x$ be a member of $\nu K$ which is
congruent to a non-trivial sum of atoms modulo $(\mu K)_r$ for
some $r\geq 0$, say
$$x\equiv\<b_1\>+\ldots+\<b_k\>\mod(\mu K)_r$$
with $k\geq 1$. If $\<b_i\>^+_r\wedge\<b_j\>^-_r=0$ for $i>j$,
then there is a decomposition
$$x=x_1\c_r\ldots\c_r x_k$$
with $x_i\in\nu K$ such that $x_i\equiv\<b_i\>\mod(\mu K)_r$.
\endproclaim

\demo{Proof} The congruence gives
$$x=\<b_1\>+\ldots+\<b_k\>+z$$
for some $z\in(\mu K)_r$; thus $d^\alpha_r z=z$. For $1\leq i\leq
k$ let $x_i$ be the element of $\mu K$ given by
$$x_i
 =d^+_r[\<b_1\>+\ldots+\<b_{i-1}\>]
 +\<b_i\>
 +d^-_r[\<b_{i+1}\>+\ldots+\<b_k\>]
 +z,$$
and for $1\leq i\leq k-1$ let
$$y_i
 =d^+_r[\<b_1\>+\ldots+\<b_i\>]
 +d^-_r[\<b_{i+1}\>+\ldots+\<b_k\>]
 +z.$$
Then $d^+_r x_i=y_i=d^-_r x_{i+1}$ because $d^\beta_r
d^\alpha_r=d^\alpha_r$, and we also have
$$x=x_1-y_1+x_2-y_2+\ldots+x_k,$$
so there is a decomposition
$$x=x_1\c_r\ldots\c_r x_k$$
in $\mu K$. It is clear that $x_i\equiv\<b_i\>\mod(\mu K)_r$, and
it remains to show that $x_i\in\nu K$. We must therefore show that
$\epsilon(x_i)^\alpha_0=1$ and that $(x_i)^\alpha_n\geq 0$.

As to the augmentation, we have
$\epsilon(d^\beta_r\<b_j\>)^\alpha_0=\epsilon\<b_j\>^\alpha_0$ for
all $j$~and~$\beta$, so $\epsilon(x_i)^\alpha_0=\epsilon
x^\alpha_0=1$.

As to the~$(x_i)^\alpha_n$, if $n<r$ then
$(x_i)^\alpha_n=x^\alpha_n\geq 0$, and if $n>r$ then
$(x_i)^\alpha_n=\<b_i\>^\alpha_n\geq 0$; it therefore remains to
consider the case $n=r$. Now $(x_1)^-_r=x^-_r\geq 0$ and
$(x_k)^+_r=x^+_r\geq 0$, and for $1\leq i\leq k-1$ we have
$$\eqalign{(x_i)^+_r
 =(x_{i+1})^-_r
 &=x^-_r+[\<b_1\>^+_r+\ldots+\<b_i\>^+_r]
 -[\<b_1\>^-_r+\ldots+\<b_i\>^-_r]\cr
 &=x^+_r+[\<b_{i+1}\>^-_r+\ldots+\<b_k\>^-_r]
 -[\<b_{i+1}\>^+_r+\ldots+\<b_k\>^+_r].\cr}$$
But $x^\alpha_r\geq 0$, $\<b_j\>^\alpha_r\geq 0$, and
$$[\<b_1\>^-_r+\ldots+\<b_i\>^-_r]
 \wedge[\<b_{i+1}\>^+_r+\ldots+\<b_k\>^+_r]
 =0,$$
so $(x_i)^\alpha_r\geq 0$ in all cases. This completes the proof.
\enddemo

To find situations in which Proposition~5.1 can be applied, we use
the following result.

\proclaim{Proposition 5.2} Let $K$ be an augmented directed
complex with a basis and let $x$ be a member of $\nu K$.

{\rm (i)} If $x\equiv 0\mod(\mu K)_{r+1}$ then
$$x\equiv\<c_1\>+\ldots+\<c_l\>\mod(\mu K)_r$$
for some $(r+1)$-dimensional atoms $\<c_1\>,\ldots,\<c_l\>$.

{\rm (ii)} If $x\equiv m\<a\>\mod(\mu K)_{r+1}$ for some positive
integer~$m$ and for some atom~$\<a\>$ with $|a|>r+1$, then
$$x\equiv m\<a\>+\<c_1\>+\ldots+\<c_l\>\mod(\mu K)_r$$
for some $(r+1)$-dimensional atoms $\<c_1\>,\ldots,\<c_l\>$.

{\rm (iii)} If $x=y\c_r z$ for some $y,z\in\nu K$ such that
$y\equiv\<a\>$ and $z\equiv\<b\>\mod(\mu K)_r$ with $|a|,|b|>r$
and $\<a\>^+_r\wedge\<b\>^-_r=0$, then
$$x\equiv\<a\>+\<b\>+\<c_1\>+\ldots+\<c_l\>\mod(\mu K)_{r-1}$$
for some $r$-dimensional atoms $\<c_1\>,\ldots,\<c_l\>$.
\endproclaim

\demo{Proof} (i) Since $x\in(\mu K)_{r+1}$, we have
$x^-_{r+1}=x^+_{r+1}=w$ for some $w\in K_{r+1}$. Since $x\in\nu
K$, we have $w\geq 0$, so that $w$~is a sum of $(r+1)$-dimensional
basis elements $c_1+\ldots+c_l$. It now follows that
$x\equiv\<c_1\>+\ldots+\<c_l\>\mod(\mu K)_r$, as required.

(ii) Here $x-m\<a\>$ is in $(\mu K)_{r+1}$. As in the proof of
part~(i), we get
$x^-_{r+1}-m\<a\>^-_{r+1}=x^+_{r+1}-m\<a\>^+_{r+1}=w$ for some
$w\in K_{r+1}$, and it suffices to show that $w\geq 0$. But this
holds because $x^-_{r+1}\geq 0$, $x^+_{r+1}\geq 0$ and
$\<a\>^-_{r+1}\wedge\<a\>^+_{r+1}=0$.

(iii) Here
$$x^-_r-\<a\>^-_r-\<b\>^-_r
 =x^+_r-\<a\>^+_r-\<b\>^+_r
 =z^-_r-\<a\>^+_r-\<b\>^-_r
 =w$$
for some $w\in K_r$, since $x\equiv y+z\equiv\<a\>+\<b\>\mod(\mu
K)_r$, since $x^+_r=z^+_r$, and since $z-\<b\>\in(\mu K)_r$. As
before, it suffices to show that $w\geq 0$. But this holds because
$x^-_r\geq 0$, $x^+_r\geq 0$, $z^-_r\geq 0$ and
$$\<a\>^-_r\wedge\<a\>^+_r
 =\<b\>^-_r\wedge\<b\>^+_r
 =\<a\>^+_r\wedge\<b\>^-_r
 =0.$$
\enddemo

To construct decompositions we use Proposition~5.2(i) together
with the case $m=1$ of Proposition~5.2(ii). It is convenient to
make the following definition.

\definition{Definition 5.3} Let $K$ be an augmented directed
complex with a basis, and let $x$ be a member of $\nu K$. Then the
{\it decomposition index\/} of~$x$ is the smallest integer $r\geq
-1$ such that $x$~is congruent to zero or an atom modulo $(\mu
K)_{r+1}$.
\enddefinition

In Definition~5.3, note that the decomposition index exists
because $x\equiv 0\mod(\mu K)_n$ for all sufficiently large~$n$.

By combining Propositions 5.1 and~5.2 we get the following result,
analogous to part of Proposition~4.2.

\proclaim{Proposition 5.4} Let $K$ be an augmented directed
complex with a loop-free basis, and let $x$ be a member of $(\nu
K)_n$ with decomposition index $r\geq 0$. Then there is a
decomposition
$$x=x_1\c_r\ldots\c_r x_k$$
with $x_i\in(\nu K)_n$ and $x_i\equiv\<b_i\>\mod(\mu K)_r$, where
$\<b_1\>,\ldots,\<b_k\>$ is a list of atoms such that $k\geq 2$,
such that $|b_i|>r$ for all~$i$, such that $|b_i|>r+1$ for at most
one value of~$i$, and such that $\<b_i\>^+_r\wedge\<b_j\>^-_r=0$
for $i>j$.
\endproclaim

\demo{Proof} We have $x$ congruent to zero or an atom modulo $(\mu
K)_{r+1}$. By Proposition~5.2(i) or~(ii) there is a congruence
$$x\equiv\<b_1\>+\ldots+\<b_k\>\mod(\mu K)_r$$
for some list of atoms~$\<b_i\>$ with $|b_i|>r$ for each~$i$ and
with $|b_i|>r+1$ for at most one value of~$i$. Since the
decomposition index is~$r$, we have $k\geq 2$. Since the basis is
loop-free we can assume the list ordered so that
$\<b_i\>^+_r\wedge\<b_j\>^-_r=0$ for $i>j$. By Proposition~5.1
there is a decomposition $x=x_1\c_r\ldots\c_r x_k$ with $x_i\in\nu
K$ and $x_i\equiv\<b_i\>\mod(\mu K)_r$. Finally, $|b_i|\leq n$ for
each~$i$, because $x\in(\nu K)_n$, and it follows that $x_i\in(\nu
K)_n$. This completes the proof.
\enddemo

To complete the analogy with Proposition~4.2 we give the following
result.

\proclaim{Proposition 5.5} Let $K$ be an augmented directed
complex with a unital basis, and let $x$ be a member of $(\nu
K)_n$ with decomposition index~$-1$. Then $x$~is an atom of
dimension at most~$n$.
\endproclaim

\demo{Proof} We have $x$ congruent to zero or an atom modulo $(\mu
K)_0$. By Proposition~5.2(i) or~(ii) there is a congruence
$$x\equiv\<b_1\>+\ldots+\<b_k\>\mod(\mu K)_{-1}$$
for some atoms~$\<b_i\>$. Since $(\mu K)_{-1}=0$, this congruence
is an equality, and it follows that
$$\epsilon x^\alpha_0
 =\epsilon\<b_1\>^\alpha_0+\ldots+\epsilon\<b_k\>^\alpha_0.$$
But $\epsilon x^\alpha_0=1$ because $x\in\nu K$, and
$\epsilon\<b_i\>^\alpha_0=1$ for each~$i$ because the basis is
unital, so $k=1$. This means that $x$~is equal to the
atom~$\<b_1\>$. The dimension of~$\<b_1\>$ is at most~$n$ because
$x\in(\nu K)_n$.
\enddemo

We can now prove one half of the equivalence.

\proclaim{Theorem 5.6} Let $K$ be an augmented directed complex
with a loop-free unital basis. Then $\nu K$ has a loop-free atomic
basis consisting of its atoms, and the counit $\pi\colon\lambda\nu
K\to K$ is an isomorphism.
\endproclaim

\demo{Proof} First we show that $\nu K$ is composition-generated
by its atoms. Indeed the atoms of dimension at most~$n$ lie in
$(\nu K)_n$ because the basis is unital, so the composites of
atoms of dimension at most~$n$ lie in $(\nu K)_n$. Conversely, let
$x$ be a member of $(\nu K)_n$ of decomposition index~$r$. If
$r\geq 0$, then $x$~is a composite of members of $(\nu K)_n$ with
lower decomposition index by Proposition~5.4, and if $r=-1$ then
$x$~is an atom of dimension at most~$n$ by Proposition~5.5. It
follows by induction on~$r$ that $x$~is a composite of atoms of
dimension at most~$n$. Therefore $\nu K$ is composition-generated
by its atoms.

Next we show that $\pi\colon\lambda\nu K\to K$ is an isomorphism.
For each $n\geq 0$ we must show that $\pi\colon(\lambda\nu K)_n\to
K_n$ is an isomorphism with $\pi[(\lambda\nu K)_n^*]=K_n^*$. Now
it follows from Proposition~4.3 that $(\lambda\nu K)_n$ and
$(\lambda\nu K)_n^*$ are generated as abelian group and submonoid
by the elements~$[\<b\>]_n$ for $|b|=n$. The images of these
elements under~$\pi$ are the $n$-dimensional basis elements, which
generate $K_n$~and~$K_n^*$ as abelian group and submonoid freely.
It follows that $\pi\colon(\lambda\nu K)_n\to K_n$ is an
isomorphism with $\pi[(\lambda\nu K)_n^*]=K_n^*$ as required.

Finally we show that the atoms form a basis for $\nu K$ with the
desired properties by verifying the conditions of Definitions 4.4
and~4.5. We have already shown that the atoms composition-generate
$\nu K$. The elements $[\<b\>]_{|b|}$ form a basis for $\lambda\nu
K$, because their images under the isomorphism
$\pi\colon\lambda\nu K\to K$ form a basis for~$K$, so the function
$\<b\>\mapsto[\<b\>]_{|b|}$ is a bijection from the set of atoms
to a basis for $\lambda\nu K$. The atoms therefore form a basis
for $\nu K$. We have $[d^-_n\<b\>]_n\wedge[d^+_n\<b\>]_n=0$ for
$n<|b|$ because
$$\pi[d^-_n\<b\>]_n\wedge\pi[d^+_n\<b\>]_n=\<b\>^-_n\wedge\<b\>^+_n=0;$$
therefore the basis for $\nu K$ is atomic. The basis for
$\lambda\nu K$ is loop-free because $\lambda\nu K\cong K$;
therefore the basis for $\nu K$ is loop-free. This completes the
proof.
\enddemo

Conversely, if $C$~is an $\omega$-category with a loop-free atomic
basis, then we want the unit $\eta\colon C\to\nu\lambda C$ to be
an isomorphism. More generally, we will prove the following
result.

\proclaim{Theorem 5.7} Let $C$ be an $\omega$-category
composition-generated by a set~$E$, let $K$ be an augmented
directed complex with a loop-free unital basis, and let
$\theta\colon C\to\nu K$ be a morphism which restricts to a
dimension-preserving bijection from~$E$ to the atoms of $\nu K$.
Then $\theta$~is an isomorphism.
\endproclaim

In a sense, Theorem~5.7 says that $\nu K$ is freely
composition-generated by its atoms. To prove it, we will show that
the decompositions in Proposition~5.4 are as nearly as possible
uniquely determined. We begin with the following observation.

\proclaim{Proposition 5.8} Let $K$ be an augmented directed
complex. If
$$x_1\c_r\ldots \c_r x_k=y_1\c_r\ldots\c_r y_k$$
in $\mu K$ and if $x_i\equiv y_i\mod(\mu K)_r$ for $1\leq i\leq
k$, then $x_i=y_i$ for $1\leq i\leq k$.
\endproclaim

\demo{Proof} We have
$$(x_1-y_1)\c_r\ldots\c_r(x_k-y_k)=0$$
and $x_1-y_1\in(\mu K)_r$, so
$$x_1-y_1=d^-_r(x_1-y_1)=d^-_r(0)=0,$$
and it follows that $x_1=y_1$. We then have
$d^-_r(x_2-y_2)=d^+_r(x_1-y_1)=0$, so $x_2=y_2$ by a similar
argument, and so on.
\enddemo

Because of Proposition~5.8, the decomposition in Proposition~5.4
is uniquely determined by the ordered list of atoms
$\<b_1\>,\ldots,\<b_k\>$. Up to permutation this list is
determined by~$x$, because $x\equiv\<b_1\>+\ldots+\<b_k\>\mod(\mu
K)_r$; however, it may be possible to reorder the list.

In particular, suppose that $k=2$, so that
$x\equiv\<b_1\>+\<b_2\>\mod(\mu K)_r$ with
$\min\{|b_1|,|b_2|\}=r+1$. Suppose further that $r>0$, and that
$x$~has a $\c_{r-1}$-decomposition $x=x_1^*\c_{r-1}x_2^*$ with
$x_1^*\equiv\<b_1\>$ and $x_2^*\equiv\<b_2\>\mod(\mu K)_r$. Then
it follows from Definition 2.1(iv) and~(vii) that there are
$\c_r$-decompositions
$$x=(x_1^*\c_{r-1}d^-_r x_2^*)\c_r(d^+_r x_1^*\c_{r-1}x_2^*)
 =(d^-_r x_1^*\c_{r-1}x_2^*)\c_r(x_1^*\c_{r-1}d^+_r x_2^*)$$
with $(x_1^*\c_{r-1}d^\alpha_r x_2^*)\equiv\<b_1\>$ and
$(d^\alpha_r x_1^*\c_{r-1}x_2^*)\equiv\<b_2\>\mod(\mu K)_r$. For
elements of this type, one can therefore transpose the atoms. We
will now show that any two consecutive atoms in the `wrong' order
can be transposed in this way.

\proclaim{Proposition 5.9} Let $K$ be an augmented directed
complex with a loop-free unital basis. Suppose that $x=y\c_r z$ is
a composite in $\nu K$ with $y\equiv\<a\>$ and $z\equiv\<b\>$
modulo $(\mu K)_r$, where $\<a\>$~and~$\<b\>$ are atoms such that
$|a|,|b|>r$ and $\<a\>^+_r\wedge\<b\>^-_r=0$. Then $r>0$ and there
is a decomposition $x=y^*\c_{r-1}z^*$ or $x=z^*\c_{r-1}y^*$ in
$\nu K$ with $y^*\equiv\<a\>$ and $z^*\equiv\<b\>$ modulo $(\mu
K)_r$.
\endproclaim

\demo{Proof} By Proposition~5.2(iii),
$$x\equiv\<a\>+\<b\>+\<c_1\>+\ldots+\<c_l\>\mod(\mu K)_{r-1}$$
for some $r$-dimensional atoms~$\<c_i\>$. If $r=0$ then this
congruence is an equality and we get $\epsilon x^\alpha_0=l+2>1$,
which is absurd. Therefore $r>0$. Since the basis is loop-free we
can order the list $\<a\>,\<b\>,\<c_1\>,\ldots,\<c_l\>$ so that
Proposition~5.1 applies. We get a decomposition
$x=x_1\c_{r-1}\ldots\c_{r-1}x_{l+2}$ in $\nu K$ such that modulo
$(\mu K)_{r-1}$ the factors are congruent to
$\<a\>,\<b\>,\<c_1\>,\ldots,\<c_l\>$ in some order. Modulo $(\mu
K)_r$ it follows that one of the factors is congruent to~$\<a\>$,
another factor is congruent to~$\<b\>$, and the others are
congruent to zero. By grouping the factors appropriately, we get a
decomposition $x=y^*\c_{r-1}z^*$ or $x=z^*\c_{r-1}y^*$ of the
required form.
\enddemo

In the proof of Theorem~5.7, we show that elements of $\nu K$ have
unique inverse images by induction on their decomposition indices.
In the inductive step we use the following lemma.

\proclaim{Lemma 5.10} Let $K$ be an augmented directed complex
with a loop-free unital basis, let $C$ be an $\omega$-category,
let $r$ be a nonnegative integer, and let $\theta\colon C\to\nu K$
be a morphism such that elements of $\nu K$ with decomposition
index less than~$r$ have unique inverse images in~$C$. If
$\xi=\eta\c_r\zeta$ is a composite in~$C$ such that
$\theta(\eta)\equiv\<a\>$ and $\theta(\zeta)\equiv\<b\>$ modulo
$(\mu K)_r$, where $\<a\>$~and~$\<b\>$ are atoms such that
$|a|,|b|>r$ and $\<a\>^+_r\wedge\<b\>^-_r=0$, then there is a
decomposition $\xi=\zeta'\c_r\eta'$ in~$C$ such that
$\theta(\zeta')\equiv\<b\>$ and $\theta(\eta')\equiv\<a\>$ modulo
$(\mu K)_r$.
\endproclaim

\demo{Proof} By Proposition~5.9, $r>0$ and there is a
decomposition $\theta(\xi)=y^*\c_{r-1}z^*$ or
$\theta(\xi)=z^*\c_{r-1}y^*$ in $\nu K$ with $y^*\equiv\<a\>$ and
$z^*\equiv\<b\>$ modulo $(\mu K)_r$. For definiteness, suppose
that $\theta(\xi)=y^*\c_{r-1}z^*$. Then $y^*$~and~$z^*$ have
decomposition indices less than~$r$, so they have unique inverse
images $\eta^*$~and~$\zeta^*$. We now get
$$\theta(d^+_{r-1}\eta^*)
 =d^+_{r-1}y^*
 =d^-_{r-1}z^*
 =\theta(d^-_{r-1}\zeta^*),$$
and the element $d^+_{r-1}y^*=d^-_{r-1}z^*$ clearly has
decomposition index less than $r-1$, so
$d^+_{r-1}\eta^*=d^-_{r-1}\zeta^*$ by the uniqueness part of the
hypothesis. It follows that $\eta^*$~and~$\zeta^*$ have a
composite $\eta^*\c_{r-1}\zeta^*$. We now have
$\theta(\xi)=\theta(\eta^*\c_{r-1}\zeta^*)$, and we claim that
$\xi=\eta^*\c_{r-1}\zeta^*$. Indeed it follows from
Definition~2.1(vii) that
$$\theta(\eta)\c_r\theta(\zeta)
 =\theta(\xi)
 =\theta(\eta^*\c_{r-1}\zeta^*)
 =\theta(\eta^*\c_{r-1}d^-_r\zeta^*)
 \c_r\theta(d^+_r\eta^*\c_{r-1}\zeta^*)$$
with
$$\theta(\eta)
 \equiv\<a\>
 \equiv y^*
 \equiv\theta(\eta^*)
 \equiv\theta(\eta^*\c_{r-1}d^-_r\zeta^*)\mod(\mu K)_r$$
and
$$\theta(\zeta)
 \equiv\<b\>
 \equiv z^*
 \equiv\theta(\zeta^*)
 \equiv\theta(d^+_r\eta^*\c_{r-1}\zeta^*)\mod(\mu K)_r,$$
so $\theta(\eta)=\theta(\eta^*\c_{r-1}d^-_r\zeta^*)$ and
$\theta(\zeta)=\theta(d^+_r\eta^*\c_{r-1}\zeta^*)$ by
Proposition~5.8. Clearly $\theta(\eta)$ and $\theta(\zeta)$ have
decomposition indices less than~$r$, so
$\eta=\eta^*\c_{r-1}d^-_r\zeta^*$ and
$\zeta=d^+_r\eta^*\c_{r-1}\zeta^*$ by the uniqueness part of the
hypothesis. Using Definition~2.1(vii) again, we get
$\xi=\eta\c_r\zeta=\eta^*\c_{r-1}\zeta^*$ as claimed. It now
follows that
$$\xi
 =(d^-_r\eta^*\c_{r-1}\zeta^*)\c_r(\eta^*\c_{r-1}d^+_r\zeta^*)
 =\zeta'\c_r\eta',$$
say, with $\theta(\zeta')\equiv\theta(\zeta^*)\equiv\<b\>$ and
$\theta(\eta')\equiv\theta(\eta^*)\equiv\<a\>$ modulo $(\mu K)_r$.
This completes the proof.
\enddemo

\demo{Proof of Theorem\/~{\rm 5.7}} We show that $\theta$~is an
isomorphism by showing that each element~$x$ of $\nu K$ has a
unique inverse image under~$\theta$. We use induction on the
decomposition index of~$x$.

Suppose first that $x$~has decomposition index~$-1$. By
Proposition~5.5, $x$~is an atom. From Proposition~4.2, we see that
the inverse images of~$x$ must be generators. Since $\theta$~maps
the set of generators bijectively onto the set of atoms, it
follows that $x$~has a unique inverse image.

Now suppose that $x$~has decomposition index $r\geq 0$. By
Proposition~5.4 there is a decomposition
$$x=x_1\c_r\ldots\c_r x_k$$
with $k\geq 2$, with $x_i\in\nu K$, and with
$x_i\equiv\<b_i\>\mod(\mu K)_r$, such that $|b_i|>r$ for all~$i$,
such that $|b_i|>r+1$ for at most one value of~$i$, and such that
$\<b_i\>^+_r\wedge\<b_j\>^-_r=0$ for $i>j$. The~$x_i$ clearly have
decomposition indices less than~$r$, so they have unique inverse
images~$\xi_i$ by the inductive hypothesis. As in the proof of
Lemma~5.10,
$$\theta(d^+_r\xi_{i-1})
 =d^+_r x_{i-1}
 =d^-_r x_i
 =\theta(d^-_r\xi_i)$$
and the decomposition index of $d^+_r x_{i-1}=d^-_r x_i$ is less
than~$r$, so $d^+_r\xi_{i-1}=d^-_r\xi_i$ by the uniqueness part of
the inductive hypothesis, and it follows that there is a composite
$$\xi=\xi_1\c_r\ldots\c_r\xi_k.$$
Then $\xi$~is an inverse image for~$x$.

To complete the proof, let $\xi'$ be any inverse image for~$x$; we
must show that $\xi'=\xi$. We can express~$\xi'$ as in
Proposition~4.2. From the form of~$x$, this must mean that
$\xi'=\xi'_1\c_r\ldots\c_r\xi'_k$ with the $\theta(\xi'_i)$
congruent to the~$\<b_i\>$ in some order modulo $(\mu K)_r$. By
repeated application of Lemma~5.10, we can change back to the
original order. This gives us a decomposition
$$\xi'=\xi''_1\c_r\ldots\c_r\xi''_k$$
with $\theta(\xi''_i)\equiv\<b_i\>\equiv x_i\mod(\mu K)_r$ for
each~$i$. Since $\theta(\xi')=x=x_1\c_r\ldots\c_r x_k$, it follows
from Proposition~5.8 that $\theta(\xi''_i)=x_i$ for each~$i$. But
the~$x_i$ have unique inverse images~$\xi_i$, so $\xi''_i=\xi_i$
for each~$i$, and it follows that $\xi'=\xi$. This completes the
proof.
\enddemo

The main theorem is now as follows.

\proclaim{Theorem 5.11} The adjoint functors
$\lambda\colon\omegacat\to\ADC$ and $\nu\colon\ADC\to\omegacat$
restrict to adjoint equivalences between the full subcategories
consisting of $\omega$-categories with loop-free atomic bases and
of augmented directed complexes with loop-free unital bases. Under
these equivalences, $\omega$-categories with strongly loop-free
atomic bases correspond to augmented directed complexes with
strongly loop-free unital bases.
\endproclaim

\demo{Proof} By Theorem~5.6, if $K$~is an augmented directed
complex with a loop-free unital basis then $\nu K$ has a loop-free
atomic basis and the counit $\pi\colon\lambda\nu K\to K$ is an
isomorphism.

Conversely, let $C$ be an $\omega$-category with a loop-free
atomic basis~$E$. By Definitions 4.4 and~4.5, the function
$e\mapsto [e]_{|e|}$ maps~$E$ bijectively to a loop-free basis for
$\lambda C$. By Proposition~4.6 this basis is unital, so $\lambda
C$ has a loop-free unital basis. By Proposition~4.6 again, the
unit $\eta\colon C\to\nu\lambda C$ sends the generator~$e$ to the
atom $\<[e]_{|e|}\>$. It follows that $\eta$~maps the generators
of~$C$ bijectively to the atoms of $\nu\lambda C$ and preserves
dimensions. By Theorem~5.7, $\eta\colon C\to\nu\lambda C$ is an
isomorphism.

Finally, strongly loop-free atomic bases correspond to strongly
loop-free unital bases by Definition 4.5(iii).

This completes the proof.
\enddemo

\section{6. Relations with earlier work}

In this section we compare our construction~$\nu$ with other
constructions of $\omega$-categories.

We first observe that Theorem~5.11 serves to characterise certain
$\omega$-categories: if $K$~is an augmented directed complex with
a loop-free basis and if $C$~is an $\omega$-category with an
atomic basis such that $\lambda C\cong K$, then $C$~must be
isomorphic to $\nu K$. This shows that the various constructions
are essentially equivalent. For example, let $C$ be the
$\omega$-category associated to the $p$-simplex in [2], [10]
or~[11]; then $C$~must be isomorphic to $\nu\Delta[p]$, where
$\Delta[p]$ is as in Example~3.8.

The $\omega$-categories associated to loop-free structures in
earlier treatments have presentations of particular types. The
analogous result for the functor~$\nu$ is as follows.

\proclaim{Theorem 6.1} Let $K$ be an augmented directed complex
with a loop-free unital basis. Then the $\omega$-category $\nu K$
has a presentation as follows\/{\rm:} the generators are the
atoms\/{\rm;} for each atom~$\<b\>$ there are relations
$d^-_{|b|}\<b\>=d^+_{|b|}\<b\>${\rm;} for each
positive-dimensional atom~$\<b\>$ there are relations
$d^-_{|b|-1}\<b\>=w^-(b)$ and $d^+_{|b|-1}\<b\>=w^+(b)$, where the
$w^\alpha(b)$ are expressions for the $d^\alpha_{|b|-1}\<b\>$ as
composites of atoms of dimension less than~$|b|$.
\endproclaim

\demo{Proof} It follows from Theorem~5.5 that $\nu K$ is
composition-generated by its atoms, and the atoms therefore
satisfy relations of the form described. Let $C$ be the
$\omega$-category generated by the atoms subject to these
relations; then there is a canonical morphism $\theta\colon
C\to\nu K$, and we must show that it is an isomorphism. Because of
Theorem~5.7, it suffices to show that the atoms
composition-generate~$C$. Because of Definition 2.1(iii) and~(v),
$C$~is composition-generated by the elements $d^\alpha_n\<b\>$ for
$\<b\>$ an atom. It therefore suffices to show that in~$C$ each
element $d^\alpha_n\<b\>$ is a composite of atoms of dimension at
most~$n$. We do this by induction on~$|b|$. Indeed, if $n<|b|$
then $d^\alpha_n\<b\>=d^\alpha_n d^\alpha_{|b|-1}\<b\>=d^\alpha_n
w^\alpha(b)$, and $w^\alpha(b)$ is a composite of atoms of
dimension less than~$|b|$, so the result holds by the inductive
hypothesis and Definition~2.1(v), and if $n\geq |b|$ then
$d^\alpha_n\<b\>=d^\alpha_n
d^\alpha_{|b|}\<b\>=d^\alpha_{|b|}\<b\>=\<b\>$, so the result
holds trivially. This completes the proof.
\enddemo

The constructions corresponding to~$\nu$ in earlier treatments are
described combinatorially in terms of sets rather than
algebraically. We will now show that $\nu K$ can be described
combinatorially when $K$~has a loop-free unital basis~$B$, in the
sense that the elements of $\nu K$ are determined by subsets
of~$B$.

\proclaim{Theorem 6.2} Let $K$ be an augmented directed complex
with a loop-free unital basis and let $x$ be a member of $\nu K$.

{\rm (i)} If $x$~is a sum of atoms then $x$~is an atom.

{\rm (ii)} If $x\equiv m_1\<b_1\>+\ldots+m_k\<b_k\>\mod(\mu K)_r$
for some $r\geq -1$, with $|b_i|>r$ and with $m_i>0$ for all~$i$,
then $m_i=1$ for all~$i$.

{\rm (iii)} Each term~$x^\alpha_n$ of~$x$ is a sum of distinct
basis elements.
\endproclaim

\demo{Proof} (i) Suppose that $x=\<b_1\>+\ldots+\<b_k\>$. Then
$\epsilon
x^\alpha_0=\epsilon\<b_1\>^\alpha_0+\ldots+\epsilon\<b_k\>^\alpha_0$.
But $\epsilon x^\alpha_0=1$ because $x\in\nu K$, and
$\epsilon\<b_i\>^\alpha_0=1$ for each~$i$ because the basis is
unital, so $k=1$. This means that $x$~is equal to the
atom~$\<b_1\>$.

(ii) The proof is by induction on~$r$. In the case $r=-1$ the
result follows from part~(i). From now on, suppose that $r\geq 0$.
We may assume that $k\geq 1$ (otherwise there is nothing to
prove). We have $\<b_i\>^+_r\wedge\<b_i\>^-_r=0$ for each~$i$,
because $|b_i|>r$. Since the basis is loop-free, we can assume the
list $\<b_1\>,\ldots,\<b_k\>$ ordered so that
$\<b_i\>^+_r\wedge\<b_j\>^-_r=0$ for $i>j$. By applying
Proposition~5.1 and grouping the factors, we get a decomposition
$x=x_1\c_r\ldots\c_r x_k$ with $x_i\in\nu K$ and $x_i\equiv
m_i\<b_i\>$ modulo $(\mu K)_r$. By Proposition~5.2(ii),
$$x_i\equiv m_i\<b_i\>+\<c_1\>+\ldots+\<c_l\>\mod(\mu K)_{r-1}$$
for some $r$-dimensional atoms $\<c_1\>,\ldots,\<c_l\>$. By the
inductive hypothesis, $m_i=1$ as required.

(iii) We have $x^\alpha_n\geq 0$ because $x\in\nu K$, so
$x^\alpha_n=m_1 b_1+\ldots+m_k b_k$ for some positive
integers~$m_i$ and some $n$-dimensional basis elements~$b_i$. It
now suffices to show that $m_i=1$ for each~$i$, and this follows
from part~(ii) because
$$d^\alpha_n x
 \equiv m_1\<b_1\>+\ldots+m_k\<b_k\>\mod(\mu K)_{n-1}.$$
\enddemo

\section{7. The structure of the category of $\omega$-categories}

The $\omega$-categories $F[p]$, $F[p;n]$, $F[p;n,n]$ and
$F[p;n,m,n]$ of Example~4.7 represent the elements, operations and
defining identities of $\omega$-categories; in particular
$F[p;n,n]$ and $F[p;n,m,n]$ are what one needs for the identities
of Definition 2.1(vi) and~(vii). One can therefore use these
$\omega$-categories to give a `globular' description of
$\omegacat$. There are also simplicial and cubical descriptions,
using the $\omega$-categories $\nu\Delta[p]$ and $\nu Q[p]$ coming
from Examples 3.8 and~3.10; see [1], [2],~[11]. All of these
descriptions are based on $\omega$-categories with strongly
loop-free atomic bases. Because of Theorem~5.11, one can express
these descriptions in terms of augmented directed complexes; thus
the theory of $\omega$-categories can be expressed in terms of
chain complexes. We will now use the globular description to get
results on monoidal structures. The main novelty is the
functoriality; see Theorem~7.6 in particular.

We first make the following observation.

\proclaim{Theorem 7.1} Every $\omega$-category is the colimit of a
small diagram of $\omega$-categories with strongly loop-free
atomic bases.
\endproclaim

\demo{Proof} Let $C$ be an $\omega$-category. By Proposition~2.3,
$C$~is the union of the increasing sequence of
sub-$\omega$-categories~$C_p$, where
$$C_p
 =d^-_p C
 =d^+_p C
 =\{\,x\in C:d^-_px=x\,\}
 =\{\,x\in C:d^+_p x=x\,\}.$$
The structure of~$C$ is given by the unary operations~$d^\alpha_n$
and the not everywhere defined binary composition
operations~$\c_n$. Now the elements and $\c_n$-composable pairs
in~$C_p$ are represented by the $\omega$-categories $F[p]$ and
$F[p;n]$ of Example~4.7, and it follows that $C$~is the colimit of
a small diagram in which the objects are copies of the $F[p]$ and
$F[p;n]$. Indeed, one needs a copy of $F[p]$ for each element
of~$C_p$; one needs a copy of $F[p;n]$ for each $\c_n$-composable
pair in~$C_p$; and one needs morphisms to determine the inclusions
$C_{p-1}\to C_p$, the operations~$d^\alpha_n$, and the composition
operations~$\c_n$. The result now holds since, by Example~4.7, the
$F[p]$ and $F[p;n]$ have strongly loop-free atomic bases.
\enddemo

From Example~3.10, the category of augmented directed complexes
with strongly loop-free unital bases has a monoidal structure
based on the usual tensor product of chain complexes. From
Theorem~5.11 we immediately get the following consequence.

\proclaim{Proposition 7.2} There is a monoidal structure on the
full subcategory of $\omegacat$ consisting of $\omega$-categories
with strongly loop-free atomic bases. It is equivalent via
$\lambda$~and~$\nu$ to the monoidal structure on the full
subcategory of $\ADC$ consisting of atomic directed complexes with
strongly loop-free unital bases.
\endproclaim

We now use Theorem~7.1 to get a monoidal structure on the category
of $\omega$-categories.

\proclaim{Theorem 7.3} There is a colimit-preserving functor
$$(C,C')\mapsto C\otimes C'\colon
 \omegacat\times\omegacat\to\omegacat$$
extending the tensor product on $\omega$-categories with strongly
loop-free atomic bases. This functor is unique up to natural
equivalence, and it determines a monoidal structure on
$\omegacat$.
\endproclaim

\demo{Proof} Because of Theorem~7.1, it makes sense to define the
functor by
$$C\otimes C'=\colim (F\otimes F'),$$
where the colimit is taken over the morphisms $F\to C$ and $F'\to
C'$ whose domains have strongly loop-free atomic bases. It is
clear that this functor preserves colimits and that it extends the
tensor product on $\omega$-categories with strongly loop-free
atomic bases. It is also clear that these properties determine the
functor up to natural equivalence. We get a monoidal structure on
$\omegacat$ because we are starting from a monoidal structure on
the subcategory.
\enddemo

We now wish to show that the monoidal structure on $\omegacat$ is
biclosed. We will use the following result.

\proclaim{Theorem 7.4} Let $\F$ be the category of
$\omega$-categories with strongly loop-free atomic bases. Then an
$\omega$-category is naturally equivalent under the functor
$C\mapsto\hom(-,C)$ to a contravariant set-valued functor on~$\F$
which takes colimits in $\omegacat$ to limits in the category of
sets.
\endproclaim

\demo{Proof} If $C$~is an $\omega$-category, then $\hom(-,C)$ is a
contravariant set-valued functor on $\omegacat$ taking colimits to
limits. It follows that the restriction $\hom(-,C)|\F$ is a
contravariant set-valued functor on~$\F$ taking colimits in
$\omegacat$ to limits in the category of sets. The axioms of
$\omega$-category theory (Definition~2.1) reduce to these
properties of $\hom(-,C)|\F$: indeed it suffices to take the
$\omega$-categories $F[p]$, $F[p;n]$, $F[p;n,n]$ and $F[p;n,m,n]$
of Example~4.7, all of which are in~$\F$. The result follows.
\enddemo

The category~$\F$ of Theorem~7.4 is equivalent under~$\lambda$ to
a subcategory of $\ADC$ by Theorem~5.11, so Theorem~7.4 yields a
description of $\omega$-categories in terms of chain complexes.
Note in particular that $\lambda$~takes colimits in $\omegacat$ to
colimits in $\ADC$ because it is a left adjoint.

From Theorem~7.4, in the standard way, we get the following
result.

\proclaim{Theorem 7.5} The monoidal structure on $\omegacat$ is
biclosed.
\endproclaim

\demo{Proof} Given $\omega$-categories $C$~and~$D$, we must find
$\omega$-categories $\HOM(C,D)$ and $\HOM'(C,D)$ such that there
are natural equivalences
$$\hom\bigl(-,\HOM(C,D)\bigr)\cong\hom(-\otimes C,D),\quad
 \hom\bigl(-,\HOM'(C,D)\bigr)\cong\hom(C\otimes-,D).$$
But it follows from Theorem~7.3 that $\hom(-\otimes C,D)$ is a
contravariant set-valued functor on $\omegacat$ taking colimits to
limits, and it then follows from Theorem~7.4 that $\hom(-\otimes
C,D)$ is naturally equivalent to $\hom\bigl(-,\HOM(C,D)\bigr)$ for
some $\omega$-category $\HOM(C,D)$. The construction of
$\HOM'(C,D)$ is similar.
\enddemo

The category of chain complexes has a well-known closed symmetric
monoidal structure, and one can check that this induces a biclosed
monoidal structure on augmented directed complexes. Indeed, let
$K$~and~$L$ be augmented directed complexes. If $n>0$ then
$$\HOM(K,L)_n=\prod_m\hom(K_m,L_{m+n}),$$
while $\HOM(K,L)_0$ consists of the pairs $(f,\epsilon f)$ such
that $f\in\prod_m\hom(K_m,L_m)$ is a chain map from~$K$ to~$L$ and
$\epsilon f$ is an integer with $(\epsilon f)(\epsilon
x)=\epsilon(fx)$ for $x\in K_0$. If
$f\in\prod_m\hom(K_m,L_{m+n})$, then we write~$f_m$ for the
component of~$f$ in $\hom(K_m,L_{m+n})$. The boundary on
$\HOM(K,L)_n$ is given by
$$(\d f)_m=\cases{
 \d\circ f_m-(-1)^n f_{m-1}\circ\d& if $m>0$,\cr
 \d\circ f_m& if $m=0$,\cr}$$
and by $\epsilon(\d f)=0$ in the case $n=1$. The augmentation on
$\HOM(K,L)_0$ is given by $(f,\epsilon f)\mapsto \epsilon f$. The
submonoid $\HOM(K,L)_n^*$ is given by the elements~$f$ of
$\prod_m\hom(K_m,L_{m+n})$ such that $f_m(K_m^*)\subset L_{m+n}^*$
for all~$m$. The definition of $\HOM'(K,L)$ is the same as for
$\HOM(K,L)$, except that the boundary formula changes to
$$(\d f)_m=\cases{
 (-1)^m(\d\circ f_m-f_{m-1}\circ\d)& if $m>0$,\cr
 \d\circ f_m& if $m=0$.\cr}$$

The biclosed monoidal structures on augmented directed complexes
and $\omega$-categories are related as follows.

\proclaim{Theorem 7.6} Let $K$~and~$L$ be augmented directed
complexes such that $K$~has a strongly loop-free unital basis.
Then there are natural isomorphisms
$$\nu\HOM(K,L)\cong\HOM(\nu K,\nu L),\quad
 \nu\HOM'(K,L)\cong\HOM'(\nu K,\nu L).$$
\endproclaim

\demo{Proof} We give the proof for the first of these
isomorphisms. Because of Theorem~7.4, it suffices to show that
there are natural bijections
$$\hom\bigl(F,\nu\HOM(K,L)\bigr)
 \cong\hom\bigl(F,\HOM(\nu K,\nu L)\bigr)$$
for $F$ an $\omega$-category with a strongly loop-free atomic
basis. But
$$\hom\bigl(F,\nu\HOM(K,L)\bigr)
 \cong\hom\bigl(\lambda F,\HOM(K,L)\bigr)
 \cong\hom(\lambda F\otimes K,L)$$
and
$$\hom\bigl(F,\HOM(\nu K,\nu L)\bigr)
 \cong\hom(F\otimes\nu K,\nu L)
 \cong\hom\bigl(\lambda(F\otimes\nu K),L\bigr)$$
because $\lambda$~is left adjoint to~$\nu$, and $\lambda F\otimes
K\cong\lambda(F\otimes\nu K)$ by Proposition~7.2. The result
follows.
\enddemo

\section{References}

[1] F. A. Al-Agl, R. Brown and R. Steiner, Multiple categories:
the equivalence of a globular and a cubical approach, Adv. Math.
170 (2002), no.~1, 71--118.

[2] F. A. Al-Agl and R. Steiner, Nerves of multiple categories,
Proc. London Math. Soc. (3) 66 (1993), no.~1, 92--128.

[3] R. Brown and P. J. Higgins, Cubical abelian groups with
connections are equivalent to chain complexes, Homology Homotopy
Appl. 5 (2003), no.~1, 49--52.

[4] S. E. Crans and R. Steiner, Presentations of omega-categories
by directed complexes, J. Austral. Math. Soc. Ser. A 63 (1997),
no.~1, 47--77.

[5] P. Gaucher, Combinatorics of branchings in higher dimensional
automata, Theory Appl. Categ. 8 (2001), no.~12, 324--376.

[6] M. Johnson, The combinatorics of $n$-categorical pasting, J.
Pure Appl. Algebra 62 (1989), no.~3, 211--225.

[7] M. M. Kapranov and V. A. Voevodsky, Combinatorial-geometric
aspects of polycategory theory: pasting schemes and higher Bruhat
orders (list of results), Cahiers Topologie G\'eom.
Diff\'erentielle Cat\'eg. 32 (1991), no.~1, 11--27.

[8] A. Patchkoria, Chain complexes of cancellative semimodules,
Bull. Georgian Acad. Sci. 162 (2000), no.~2, 206--208.

[9] A. J. Power, An $n$-categorical pasting theorem, in Category
theory (Como, 1990), 326--358, Lecture Notes in Math. 1488
Springer, Berlin, 1991.

[10] R. Steiner, The algebra of directed complexes, Appl. Categ.
Structures 1 (1993), no.~3, 247--284.

[11] R. Street, The algebra of oriented simplexes, J. Pure Appl.
Algebra 49 (1987), no.~3, 283--335.

[12] R. Street, Parity complexes, Cahiers Topologie G\'eom.
Diff\'erentielle Cat\'eg. 32 (1991), no.~4, 315--343; Corrigenda,
Cahiers Topologie G\'eom. Diff\'erentielle Cat\'eg. 35 (1994),
no.~4, 359--361.

\bye